\documentclass[12pt]{amsart}
\usepackage{enumerate}
\usepackage{amssymb}
\usepackage{bbm}

\newenvironment{excise}[1]{}{}

\newtheorem{thm}{Theorem}[section]
\newtheorem{lemma}[thm]{Lemma}
\newtheorem{prop}[thm]{Proposition}
\newtheorem{cor}[thm]{Corollary}
\newtheorem{conj}[thm]{Conjecture}

\theoremstyle{definition}
\newtheorem{defn}[thm]{Definition}
\newtheorem{example}[thm]{Example}
\newtheorem{remark}[thm]{Remark}

\newcommand\0{\mathbf{0}}
\newcommand\1{\mathbbm{1}}
\newcommand\NN{\mathbb{N}}
\newcommand\RR{\mathbb{R}}
\newcommand\ZZ{\mathbb{Z}}
\newcommand\cB{\mathcal{B}}
\newcommand\cD{\mathcal{D}}
\newcommand\cN{\mathcal{N}}
\newcommand\cP{\mathcal{P}}
\newcommand\ve{\varepsilon}
\newcommand\ttt{\mathbf{t}}

\newcommand\minus{\smallsetminus}
\newcommand\nothing{\varnothing}
\renewcommand\iff{\Leftrightarrow}
\renewcommand\implies{\Rightarrow}

\newcommand\ceil[1]{\lceil{#1}\rceil}

\begin{document}

\mbox{}\vspace{-2.5ex}
\title{Lattice point methods for combinatorial games}
\author{Alan Guo}
\address{Department of Mathematics\\Duke University\\Durham, NC 27708}
\email{axg@duke.edu}
\author{Ezra Miller}
\address{Department of Mathematics\\Duke University\\Durham, NC 27708}
\email{ezra@math.duke.edu}

\begin{abstract}
We encode arbitrary finite impartial combinatorial games in terms of
lattice points in rational convex polyhedra.  Encodings provided by
these \emph{lattice games} can be made particularly efficient for
octal games, which we generalize to \emph{squarefree games}.  These
additionally encompass all heap games in a natural setting, in which
the Sprague--Grundy theorem for normal play manifests itself
geometrically.  We provide an algorithm to compute normal play
strategies.

The setting of lattice games naturally allows for mis\`ere play, where
$0$ is declared a losing position.  Lattice games also allow
situations where larger finite sets of positions are declared losing.
Generating functions for sets of winning positions provide data
structures for strategies of lattice games.  We conjecture that every
lattice game has a \emph{rational strategy}: a rational generating
function for its winning positions.  Additionally, we conjecture that
every lattice game has an \emph{affine stratifi\-cation}: a partition
of its set of winning positions into a finite disjoint union of
finitely generated modules for affine semigroups.
\end{abstract}

\subjclass{(MSC date is wrong: these are MSC 2010).  Primary: 91A46,
91A05, 52B20, 05A15; Secondary: 05E40, 06F05, 20M14, 68W30}



\date{3 August 2009}

\maketitle

\section{Introduction}

Combinatorial games involve two players, both with complete
information, taking turns moving on a fixed game tree.  The games
considered here are impartial, meaning that both players have the same
available moves from each position, and finite, meaning that the game
tree is finite, although we are interested in families of games in
which the totality of the available positions is infinite.  The
quintessential example of such (a family of) games is \textsc{Nim}, in
which each node of the tree corresponds to a finite set of heaps of
given sizes, and a move is accomplished by removing any number of
beans from a single heap.

The \emph{normal play} version of \textsc{Nim}, where the last person
to move is the winner, was solved over a century ago \cite{nim}.  A
complete structure theory for normal play games, known as the
Sprague--Grundy theorem from the late 1930s \cite{Sprague36,Grundy39},
builds on Bouton's solution by reducing all finite impartial games
to~it: every impartial game under normal play is equivalent to a
single \textsc{Nim} heap of some size.  Background and details can be
found in \cite{lessons,winningWays1}.

In contrast, in \emph{mis\`ere play} the last player to move is the
loser, as in Dawson's chess \cite{Dawson34}.  Mis\`ere games are much
more complex than normal play games, as observed by Conway (see
\cite{ONAG}), after inherent difficulties effectively stymied all
progress.  In particular, Dawson's chess remains elusive, despite
exciting recent advances in mis\`ere theory by Plambeck and Siegel
\cite{Pla05, misereQuots}, about which we say more shortly.

Our goal is to provide a setting in which algorithmic and theoretical
techniques concerning lattice points in polyhedra can be brought to
bear on computational and abstract periodicity questions from mis\`ere
combinatorial game theory.  Our approach is to rephrase the language
of heap-based games, especially the historically popular \emph{octal
games} (of which Dawson's chess is an example), with the consequence
of placing them in a certain natural generality, that of
\emph{squarefree games} (see Definition~\ref{d:squarefree}, and
Section~\ref{s:squarefree} in general).  The lattice encoding allows
for both normal and mis\`ere play, as well as generalizations in which
finitely many positions are declared to be losing positions
(Definition~\ref{d:board}).  The lattice encoding of squarefree games
is also efficient; that said, it turns out that our lattice point
language captures arbitrary finite impartial combinatorial games
(Theorem~\ref{t:arbitrary}), although the algorithmic
efficiency---which is key to the goal of solving octal games---is lost
in translation.

Generally speaking, any useful notion of a ``solution'' or
``strategy'' for a game should be a data structure with at least two
fundamental properties:\vfill
\begin{itemize}
\item%
it can be efficiently stored, and
\item%
it can be efficiently processed to compute a winning move from any
position.
\end{itemize}\vfill
Ideally, the data structure should also be efficiently computable,
given the game board and rule set, but that confronts a separate issue
of finding a strategy, as opposed to simply recognizing one.

In the above sense, the \emph{mis\`ere quotients} invented and mined
by Plambeck and Siegel \cite{Pla05, misereQuots} constitute excellent
data structures for solutions of mis\`ere play impartial combinatorial
games (and impartial ones \cite{allenThesis}, too), when the quotients
are finite.  This has been the basis for substantial advances in
computations involving specific mis\`ere games, such as some
previously unyielding octal games, including some with infinite
quotients; see \cite{misereQuots, Wei09}, for example.  However, even
when an infinite mis\`ere quotient is given as a finitely presented
monoid, say by generators and relations, there remains a need to
record the winning and losing positions---that is, the
\emph{bipartition} of the mis\`ere quotient into P-positions and
N-positions.

The \emph{lattice games} that we develop are played on \emph{game
boards} constituting sets of lattice points in polyhedra
(Definition~\ref{d:board}).  For example, when the polyhedron is a
cone, the game board is an \emph{affine semigroup}.  The allowed moves
constitute a finite set of vectors called a \emph{rule set}
(Definition~\ref{d:ruleset}), satisfying some conditions of
compatibility with the game board.  A~rule set uniquely determines the
sets of winning and losing positions on a given game board
(Theorem~\ref{t:uniqueness}).

Methods from combinatorial commutative algebra, as it relates to sets
of lattice points, provide clues as to how to express the presence of
structure in the sets of winning positions of games.  We conjecture
that every lattice game has
\begin{itemize}
\item%
a \emph{rational strategy}: the generating function for its winning
positions is a ratio of polynomials with integer coefficients
(Conjecture~\ref{conj:rational}); and
\item%
an \emph{affine stratification}: an expression of its winning
positions as a finite disjoint union of finitely generated modules for
affine semigroups (Conjecture~\ref{conj:strat}).
\end{itemize}
The second conjecture is stronger, in that it implies the first in an
efficient algorithmic sense (Theorem~\ref{t:strat}), but it is the
first that posits a data structure exhibiting the two fundamental
properties listed above for a successful strategy
(Theorem~\ref{t:rational}).  That said, affine stratifications reflect
real, observed phenomena more subtle than mere rationality of
generating functions.  We intend affine stratifications to provide
vehicles for producing rational strategies, and they have already been
useful as such in examples.

The above conjectures are true for squarefree games under normal play
(Corollary~\ref{c:normalplay} and Example~\ref{e:normalplay}).
Additionally, for squarefree normal play games, we present an
algorithm for computing the winning positions (Theorem~\ref{t:alg}).

\subsection*{Acknowledgements}

The authors are grateful to Mike Weimerskirch, who participated in
early discussions on this material, leading to the preliminary
definitions and conjecture in~\cite{gmw}.  A hearty thanks also goes
to Dalhousie University and all of the participants of Games-At-Dal-6
(July, 2009), organized and graciously hosted by Richard Nowakowski,
for many crucial extended discussions.  The seeds for this work were
sown at the Banff International Research Station (BIRS) Combinatorial
Game Theory Workshop, in January 2008.  AG was funded as a Duke
University undergraduate PRUV fellow for the duration of this work.
EM was partially funded by NSF Career Grant DMS-0449102.

We feel it is most appropriate for this paper to appear in a volume
dedicated to Dennis Stanton.  In fact, the mere thought of submitting
this work to his volume benefited us greatly: it brought to mind
Dennis's work on the combinatorics of various kinds of generating
functions, which led us quickly to rational functions as finite data
structures for solutions of combinatorial games.  This path took us
straight to the formulation of Conjecture~\ref{conj:rational} as a
consequence of the affine stratifications in Definition~\ref{d:strat}
and Conjecture~\ref{conj:strat}.

\section{Lattice games}\label{s:lattice}

General abstract lattice games are played on game boards consisting of
lattice points in polyhedra.  Readers unfamiliar with polyhedra can
find an elementary exposition of the necessary facts in Ziegler's book
\cite[Lectures~0--2]{ziegler}.  A \emph{polytope} is the convex hull
of a finite set in $\RR^d$.  A \emph{polyhedron} is the intersection
of finitely many closed affine halfspaces, each one bounded by a
hyperplane that need not pass through the origin. Every bounded
polyhedron is a polytope; this fact is intuitively true but not
trivial to prove \cite[Theorem~1.1]{ziegler}.

Much of the framework here will depend on the geometry and algebra of
affine semigroups for which a general reference is
\cite[Chapter~7]{cca}.  For now we recall some basic definitions.  A
\emph{semigroup} is a set with an associative binary operation.  If
the operation has an identity element, then the semigroup is a
\emph{monoid}.  An \emph{affine semigroup} is a monoid that is
isomorphic to a finitely generated subsemigroup of $\ZZ^d$ for some
$d$.  An affine semigroup is \emph{pointed}\/ if the identity is its
only unit (i.e., invertible~element).

Fix a pointed rational convex polyhedron $\Pi \subset \RR^d$ with
recession cone~$C$ of dimension~$d$.  The \emph{pointed} condition
means that $\Pi$ possesses a vertex; equivalently, $C$ is a pointed
cone.  The \emph{rational} condition means that the linear
inequalities defining $\Pi$ have rational coefficients.  Write
$\Lambda = \Pi \cap \ZZ^d$ for the set of integer points in~$\Pi$.

\begin{example}\label{e:nd}
The case of primary interest is $\Pi = C = \RR_+^d$, so $\Lambda =
\NN^d$, in which lattice points with nonnegative coordinates represent
positions in the game.  The class of \emph{heap games} is subsumed in
this context: from an initial finite set of heaps of beans, the
players take turns changing a heap---whichever they select---into some
number of heaps of smaller sizes.  The rules of a heap game specify
the allowed smaller sizes.  The game of \textsc{Nim} follows this
pattern: a player is allowed to remove beans from any single heap,
thus either creating one heap of strictly smaller size or deleting the
heap entirely.  In terms of lattice games, a position $p =
(\pi_1,\ldots,\pi_d) \in \NN^d$ represents $\pi_i$ heaps of size~$i$
for $i = 1,\ldots,d$.  Octal games, quaternary games, hexadecimal
games, and so on, are heap games; we will examine these later (under
normal play) in the wider context of \emph{squarefree games}, to be
defined and analyzed in Section~\ref{s:squarefree}.
\end{example}

Moves in lattice games will require some hypotheses in order to
guarantee that positions can reach a suitable neighborhood of the zero
game.  The geometric condition implying this behavior involves the
following notion.

\begin{defn}\label{d:tangent}
Given an extremal ray $\rho$ of a cone~$C$, the \emph{negative tangent
cone} of $C$ along~$\rho$ is $-T_\rho C = -\bigcap_{H \supset \rho}
H_+$ where $H_+ \supseteq C$ is the positive closed halfspace bounded
by a supporting hyperplane $H$ for $C$.
\end{defn}

Throughout this paper, set $\0 = (0,\ldots,0) \in \ZZ^d$.

\begin{defn}\label{d:ruleset}
A finite subset $\Gamma \subset \ZZ^d \minus \{\0\}$ is a \emph{rule
set} if
\begin{enumerate}
\item%
there exists a linear function on~$\RR^d$ that is positive on $C
\minus \{\0\}$ and on~$\Gamma$; and
\item%
for each ray $\rho$ of $C$, some vector $\gamma_\rho \in \Gamma$ lies
in the negative tangent cone $-T_\rho C$.
\end{enumerate}
\end{defn}

\begin{example}
With notation as in Example~\ref{e:nd}, the positions of the game
\textsc{Nim} with heaps of size at most~$2$ correspond to~$\NN^2$.
Each move either removes a $1$-heap, removes a $2$-heap, or turns a
$2$-heap into a $1$-heap.  Hence the rule set $\Gamma$ consists of
$(1,0)$, $(0,1)$, and $(-1,1)$, respectively.  The options of $p =
(\pi_1,\pi_2)$ are the elements of the set $(p - \Gamma) \cap \NN^2$.
We verify that $\Gamma$ is a rule set: for condition~1, the function
$\ell: \ZZ^2 \rightarrow \ZZ$ defined by $\ell(x,y) = x + 2y$ is
positive on $\NN^2 \minus \{\0\}$ and on~$\Gamma$; condition~2 is
satisfied by the basis vectors in~$\Gamma$.
\end{example}

\begin{example}[Heap games]\label{e:heaps}
In the situation of Example~\ref{e:nd}, the rule set of a heap game
is, by definition, a finite set of vectors~$\gamma$ each having the
property that all of the nonzero entries of~$\gamma$ are negative,
except the last nonzero entry of~$\gamma$, which equals~$1$.
Therefore, any linear function $\ell = (\ell_1,\ldots,\ell_d)$ is
positive on $\NN^d \minus \{\0\}$ and on~$\Gamma$ as long as $\ell_i$
is positive and sufficiently bigger than $\ell_{i-1}$ for each~$i$.
The tangent cone axiom is satisfied by definition for heaps of a given
size~$i$ as long as there is a way to act on a heap of that size; that
is, as long as some $\gamma \in \Gamma$ has $\gamma_i = 1$.
\end{example}

\begin{lemma}\label{l:pointed}
The affine semigroup $\NN\Gamma$ generated by any rule set $\Gamma$ is
pointed.
\end{lemma}
\begin{proof}
The nonzero vectors in $\RR_+ \Gamma$ all lie on the same side of the
hyperplane given by the vanishing of the linear function.
\end{proof}

\begin{remark}
Condition 1 in Definition~\ref{d:ruleset} implies more than
Lemma~\ref{l:pointed}: it implies also that $\NN\Gamma$ and $\Lambda$
point in the same direction.  That is, moves, which are elements of
$-\Gamma$, bring positions closer to~$\0$.
\end{remark}

\begin{lemma}\label{l:gammaorder}
Any rule set $\Gamma$ induces a partial order $\preceq$ on~$\Lambda$
with $p \preceq q$ if $q - p \in \NN\Gamma$.
\end{lemma}
\begin{proof}
This follows immediately from the definitions of poset and
Lemma~\ref{l:pointed}.
\end{proof}

\begin{defn}[Lattice games]\label{d:board}
Given the polyhedral set $\Lambda = \Pi \cap \ZZ^d$, fix a
rule~set~$\Gamma$.
\begin{itemize}
\item%
A \emph{game board} $\cB$ is the complement in $\Lambda$ of a finite
$\Gamma$-order ideal in $\Lambda$ called the set of \emph{defeated
positions}.
\item%
A \emph{lattice game} is defined by a game board and a rule set.
\item%
A position $p \in \Lambda$ \emph{has a move to} $q \in \Lambda$ if $p
- q \in \Gamma$.
\item%
A move from a position $p$ to~$q$ is \emph{legal} if $q$ lies on the
game board~$\cB$.
\item%
The \emph{options} of a position are the positions to which it has
legal moves.
\end{itemize}
\end{defn}

\section{Geometry of rule sets}\label{s:geometry}

The axioms for rule sets in Definition~\ref{d:ruleset} and game boards
in Definition~\ref{d:board} were chosen to define what we suspect is
the widest class of games satisfying the conjectures in
Section~\ref{s:rational}, among games played on lattice points in this
manner.  This particular choice of axioms, however, resulted from
numerous discussions about the properties of more inclusive and more
restrictive classes of games.  The purpose of this section is to
explain our rationale, including consequences of the axioms in
Section~\ref{s:lattice}, while introducing some potentially useful
extra conditions to put on rule sets.

\subsection{Why should the rule set be pointed?}\label{sub:pointed}

If rule sets were allowed to generate cones with nontrivial lineality,
then any position far from the game board boundary would have a loop
of moves back to itself, violating the blanket finiteness condition.
(This is the same reason $\0$ is not allowed as a move.)  Indeed,
since a subset of~$\Gamma$ generates the lineality group as a monoid,
some positive combination of moves along the lineality equals~$\0$.
Thus nontrivial lineality is ruled out, so we must require the cone
$\RR_+\Gamma$ to be pointed in Definition~\ref{d:ruleset}.
Consequently, just as abstract finite combinatorial games allow
induction on options, the pointed hypothesis allows for induction on
options in lattice games, by choosing a linear function~$\ell$ on the
game board that is positive on~$\Gamma$, since $\ell(p) <
\ell(p-\gamma)$ for any position~$p$.

\subsection{Why arbitrary polyhedra?}\label{sub:polyhedral}

Given the interpretation of $\NN^d$ in terms of heaps, it seemed
natural, at first, to play all lattice games on~$\NN^d$, after
disallowing $\0$ for the purposes of mis\`ere play.  However, we saw
no a~priori reason to require that a rule set must span~$\RR^d$.  But
without such a ``fullness'' hypothesis on~$\Gamma$, the positions
reachable from a given initial position~$p$ are restricted to lie in
the coset through~$p$ of the lattice $\ZZ\Gamma$ generated
by~$\Gamma$.  Thus we were led to arbitrary polyhedral game boards,
because $(p + \ZZ\Gamma) \cap \NN^d$ comprises the $\ZZ\Gamma$-lattice
points in the polyhedron $\Pi = (p + \RR\Gamma) \cap \RR_+^d$.
Although all of the examples that currently interest us are played on
the polyhedron~$\RR_+^d$, we have no reason to believe that the
conjectures in Section~\ref{s:rational} fail for arbitrary polyhedral
game boards.

\subsection{Why the tangent cone axiom?}\label{sub:tangent}

Suppose, for the moment, that we are to play a game on the
board~$\NN^d$.  If the cone $\RR_+\Gamma$ generated by the rule
set~$\Gamma$ fails to contain the nonnegative orthant~$\RR_+^d$, then
the positions in $\RR_+^d \minus \RR_+\Gamma$ can certainly never
reach the winning position~$\0$ by any sequence of moves.  Thus, when
the polyhedron $\Pi$ is the nonnegative orthant~$\RR_+^d$, we
considered using the condition $\RR_+\Gamma \supseteq \RR_+^d$ as a
rule set axiom.  The generalization of this requirement to arbitrary
polyhedral game boards is $\RR_+\Gamma \supseteq C$: the rule set cone
must contain the recession cone of the game board.  Indeed, for
general polyhedral boards, writing $\Pi = \Pi_0 + C$ for a
polytope~$\Pi_0$,~we~get
$$%
  \Pi_0 + \RR_+\Gamma \supseteq \Pi_0 + C = \Pi,
$$
so the bounded convex set $\Pi_0$ occupies the role for arbitrary
polyhedral boards that the game~$\0$ occupies for~$\NN^d$.  Note that
it is unreasonable to expect a sequence of moves to reach a single
goal position on an arbitrary polyhedral board (which position would
it be?), and this becomes especially true in generalized mis\`ere
play, where a finite subset of the lattice points are declared to be
defeated positions.

As natural and simple as it may be, the condition $\RR_+\Gamma
\supseteq C$ guarantees only that every position has a sequence of
moves to near the boundary of the game board, not necessarily to the
origin or anywhere near it (or, for arbitrary polyhedra,
near~$\Pi_0$).

\begin{example}\label{e:path}
The set $\Gamma = \{(1,0,0), (0,1,0), (1,-1,1), (-1,1,1)\}$ generates
a real pointed cone containing the nonnegative octant~$\RR_+^3$.
If~$\Gamma$ were allowed as a rule set, then no position along the
third axis would have any legal moves.
\end{example}

The example demonstrates an elementary observation: if every position
is to have a sequence of legal moves that ends in a neighborhood of
the origin, then every position sufficiently far along every extreme
ray of the game board must have at least one legal move.  The tangent
cone axiom is precisely what guarantees this; and once it does, every
position has a sequence of moves to a bounded neighborhood.  Let us be
precise.

\begin{defn}\label{d:path}
Given a finite set $\Gamma \subset \ZZ^d \minus \{\0\}$, a
\emph{$\Gamma$-path} is a sequence $p_0, \ldots, p_r$ with $p_{k-1} -
p_k \in \Gamma$ for $k = \{1,\ldots,r\}$.  For any set $S \subseteq
\ZZ^d$ of lattice points (such as a game board or the lattice points
in a polyhedron), this $\Gamma$-path \emph{goes from $p$ to~$q$
in~$S$} if $p = p_0$ and $p_r = q$ and $p_0,p_1,\ldots,p_r$ all lie
in~$S$.
\end{defn}

\begin{thm}\label{t:path}
In any lattice game, there is a finite set $F \subset \cB$ of game
board positions such that every position in~$\cB$ has a $\Gamma$-path
in~$\cB$ to~$F$.  Equivalently, the set of \emph{victorious positions}
(those from which there are no legal moves) is finite.
\end{thm}
\begin{proof}
The equivalence is straightforward, so we prove only the second claim.

There is a finite set $\Lambda_0 \subset \Lambda$ of lattice points
such that $\Lambda = \Lambda_0 + (C \cap \ZZ^d)$; indeed, one can take
for $\Lambda_0$ the set of all lattice points in a sufficiently large
neighborhood (in~$\Pi$) of the polytope~$\Pi_0$ from
Section~\ref{sub:tangent}.  Since the legal moves in~$\cB$ between
positions within $\lambda + (C \cap \ZZ^d)$ are exactly the legal
moves in~$\lambda + (C \cap \ZZ^d)$, it suffices to treat the case
$\Pi = \lambda + C$.  After translating by $-\lambda$ if necessary, we
may assume~that~$\Pi = C$.

There are two types of victorious positions in~$C$: those from which
all moves to positions in~$C$ land in~$\cD$, and those from which no
moves land in~$C$.  As the first of these sets is finite, because
$\cD$ and~$\Gamma$ are finite, we can and do assume that $\cD =
\nothing$.

Let $R$ be the set of extremal rays of~$C$.  It is enough to show that
$\bigcup_{\rho \in R} (\gamma_\rho + C)$ covers all but a bounded
subset of~$C$, because every position in $\gamma_\rho + C$ has the
legal move~$\gamma_\rho$.  Write $C_\rho$ for the union of the facets
of~$C$ that do not contain~$\rho$, and let $\ve$ be the maximum of the
lengths of the vectors~$\gamma_\rho$.  Then $\gamma_\rho + C$ contains
(perhaps properly) the set of all points of~$C$ lying outside of the
$\ve$-neighborhood of~$C_\rho$.

\begin{lemma}
Let $X$ be a polytope.  For each vertex\/~$v$ of~$X$, let $X_v$ be the
union of the facets of~$X$ not meeting~$v$.  For $\mu \in \RR_+$, let
$N_v(\mu) \subset \mu X$ be the open subset consisting of all points
that do not lie within~$\ve$ of~$\mu X_v$.  Then $\mu X = \bigcup_v
N_v(\mu)$ for all $\mu \gg 0$.
\end{lemma}
\begin{proof}
In the barycentric subdivision of~$X$, if $U_v$ denotes the union of
the closed simplices containing the vertex~$v$, then $X = \bigcup_v
U_v$.  Choose $\mu$ big enough so that the barycenter of each
positive-dimensional face~$Y$ of~$\mu X$ lies at distance greater
than~$\ve$ from the affine spans of all of the facets of~$Y$.  Then
$\mu U_v$ is contained in~$N_v(\mu)$, so $\bigcup_v N_v(\mu) \supseteq
\bigcup_v \mu U_v = X$, completing the proof of the lemma.
\end{proof}

\noindent
Using the lemma on any transverse hyperplane section~$X$ of~$C$ proves
the theorem.
\end{proof}

We stress that the conclusion of Theorem~\ref{t:path} is equivalent to
the tangent cone axiom in Definition~\ref{d:ruleset}, because of the
elementary observation following Example~\ref{e:path}.  Since many
rule sets naturally satisfy this alternate condition (or a stronger
one: the heap rule sets in Example~\ref{e:heaps} always result
in~$\Gamma$-paths to~$\0$), we record the remaining implication for
future reference.

\begin{lemma}\label{l:path}
Suppose $\Gamma$ satisfies the positivity axiom of
Definition~\ref{d:ruleset}, but not necessarily the tangent cone
axiom.  If there is finite set $F \subset \Lambda$ such that every
position in~$\Lambda$ has a $\Gamma$-path in~$\Lambda$ to~$F$, then
$\Gamma$ is a rule set (i.e., $\Gamma$ satisfies the tangent
cone~axiom).
\end{lemma}

Infinitude of the set of victorious positions (prevented by
Theorem~\ref{t:path}) could have the potential to break
Conjecture~\ref{conj:rational} on the existence of rational
strategies.  That is one of our key reasons for using the tangent cone
axiom instead of the (seemingly more) natural cone containment
condition.  With that in mind, let us prove that the tangent cone
axiom is indeed stronger.

\begin{prop}\label{p:outerTangent}
The cone generated by any rule set contains the recession cone of the
game board: $\RR_+\Gamma \supseteq C$.
\end{prop}
\begin{proof}
If $C'$ is a cone containing~$C$, and if $\rho$ is an extremal ray
of~$C$ that remains extremal in~$C'$, then automatically $-T_\rho C'
\supseteq -T_\rho C$.  It is therefore enough to show that adding a
new generating ray to~$C$ or replacing an extremal ray~$\rho$ of~$C$
by a ray~$\gamma_\rho$ in~$-T_\rho C$ yields a cone~$C'$
containing~$C$.  This is obvious when a new generating ray is added,
or when $\gamma_\rho$ lies along~$\rho$.  In the other case, the
segment from $\gamma_\rho$ through any point on~$\rho$ extends to pass
through a point~$\beta_\rho$ in some boundary face of~$C$ not
meeting~$\rho$.  By construction, $\gamma_\rho$ plus some (uniquely
defined) multiple of $\beta_\rho$ lies along~$\rho$; the positivity
axiom guarantees that it lies along $\rho$ and not~$-\rho$.  Since the
vector $\beta_\rho$ is a positive combination of extremal rays of~$C'$
none of which is~$\rho$, the ray $\rho$ remains in the new cone~$C'$
generated by $\gamma_\rho$ along with the other rays
of~$C$,~whence~\mbox{$C' \supseteq C$}.
\end{proof}

\subsection{The index of the rule set lattice}\label{sub:index}

In the general polyhedral setting, if we are interested in a fixed
starting position, then it imposes no restriction to assume that
$\Gamma$ is saturated, in the following sense, for otherwise we may
simply choose a smaller lattice to call~$\ZZ^d$.

\begin{defn}
A rule set~$\Gamma$ is \emph{saturated} if it spans~$\ZZ^d$ as a
group: $\ZZ\Gamma = \ZZ^d$.
\end{defn}

In particular, although the index of $\ZZ\Gamma$ in~$\ZZ^d$
contributes to the computational complexity, it has no effect on the
rationality or stratification conjectures in Section~\ref{s:rational},
since any lattice game breaks up into a disjoint union of
$|\ZZ^d/\ZZ\Gamma|$ many lattice games.

That said, many natural lattice games on~$\NN^d$---where replacing the
polyhedron or the lattice with new ones is undesirable---have
saturated rule sets.  Such is the case with the heap games in
Example~\ref{e:heaps}, for instance, by the following general
criterion.

\begin{prop}\label{p:saturated}
Fix a lattice game with game board~$\NN^d$ and rule set~$\Gamma$.  If
every position in~$\NN^d$ has a $\Gamma$-path to~$\0$ in~$\NN^d$, then
$\Gamma$ is saturated.
\end{prop}
\begin{proof}
The hypothesis implies that all points in~$\NN^d$ lie in the same
coset of~$\ZZ\Gamma$ as~$\0$.
\end{proof}

\section{Uniqueness of winning and losing positions}

In this section, we show that the sets of winning and losing positions
of a lattice game are well-defined.  This result, and all of the
others in this section, hold in full without the tangent cone axiom
for rule sets in Definition~\ref{d:ruleset}.  To begin, here is an
algebraic---and seemingly non-recursive---definition of winning and
losing~positions.

\begin{defn}\label{d:NP}
If $G$ is a lattice game with game board~$\cB$
and rule set~$\Gamma$, then $\cP$ is the set of \emph{winning
positions} of~$G$, and $\cN$ is the set of \emph{losing positions}
of~$G$, if $\cP$ and~$\cN$ partition $\cB$ and $(\cP + \Gamma) \cap
\cB = \cN$.
\end{defn}

In game-theoretic terms, the equation $(\cP + \Gamma) \cap \cB = \cN$
says that every position on the game board with a move to a winning
position is a losing position.  This implies other game-theoretic
statements about winning and losing positions, such as the following,
which looks similar but is strictly weaker.

\begin{prop}\label{p:p-G}
If $\cB$ is a game board with winning positions $\cP$, losing
positions~$\cN$, and rule set $\Gamma$, then $(\cP - \Gamma) \cap \cP
= \nothing$.
\end{prop}
\begin{proof}
Suppose $p \in (\cP - \Gamma) \cap \cP$.  Then $p = p' - \gamma$ for
some $p' \in \cP$ and some $\gamma \in \Gamma$.  Therefore $p' = p +
\gamma \in (\cP + \Gamma) \cap \cB = \cN$, a contradiction.
\end{proof}

Proposition~\ref{p:p-G} is weaker than Definition~\ref{d:NP} because
it does not force each losing position to possess a move to some
winning position.  For example, it is possible to change all but
finitely many P-positions to N-positions without violating
Proposition~\ref{p:p-G}, but Definition~\ref{d:NP} guarantees the
existence of infinitely many P-positions.

Next we explore consequences of the compatibility of the rule set and
the game board dictated by the positivity axiom in
Definition~\ref{d:ruleset}.

\begin{lemma}\label{l:dcc}
Every sequence in $\Lambda$ decreasing with respect to a rule set
$\Gamma$ is finite.
\end{lemma}
\begin{proof}
Fix an integer linear function $\ell$ on~$\ZZ^d$ that is positive
on~$\Gamma$.  Then $\ell(q) -\nolinebreak\ell(p) = \ell(q - p) \geq 0$
whenever $p \preceq q$, with equality if and only if $p = q$.
Therefore $\ell$ induces a bijection from each $\Gamma$-decreasing
sequence in $\Lambda$ to some decreasing sequence in~$\ZZ$.
Therefore, it is enough to show that $\ell(\Lambda)$ is bounded below.
This can be accomplished using Proposition~\ref{p:outerTangent} (or
simply by assuming that $\ell$ is nonnegative on~$C$, using
Definition~\ref{d:ruleset}.1, if one wishes to avoid invoking the
tangent cone axiom).
\end{proof}

\begin{defn}
Let $T \subseteq \ZZ^d$.  An element $p \in T$ is
\emph{$\Gamma$-minimal} if $(p - \Gamma) \cap T = \nothing$.  An
element $p \in T$ is \emph{$\NN\Gamma$-minimal} (or simply
\emph{minimal}) if $(p - \NN\Gamma) \cap T = \nothing$.
\end{defn}

\begin{example}
In mis\`ere play, we assume $\Pi = C = \RR_+^d$, $\Lambda = \NN^d$,
and $\cD = \{\0\}$.  In this case, every position has a $\Gamma$-path
to~$\0$ in~$\Lambda$ if and only if every $\Gamma$-minimal element of
$\cB$ lies in $\Gamma$.  Indeed, if a $\Gamma$-minimal element $p \in
\cB$ does not lie in $\Gamma$, then $p - \gamma \notin \NN^d$ for
every $\gamma \in \Gamma$, and hence $p$ does not have a $\Gamma$-path
to~$\0$.  Conversely, suppose every $\Gamma$-minimal element of $\cB$
lies in $\Gamma$.  By Lemma~\ref{l:dcc}, every $p \in \cB$ has a
$\Gamma$-path in~$\NN^d$ to a $\Gamma$-minimal element, and hence to
$\0$.
\end{example}

\begin{thm}\label{t:uniqueness}
Given a rule set $\Gamma \subset \mathbb{Z}^d$ and a game board~$\cB$,
there exist unique sets $\cP$ and $\cN$ of winning and losing
positions~for~$\cB$.
\end{thm}
\begin{proof}
By Lemma~\ref{l:gammaorder} and Lemma~\ref{l:dcc}, $\cB$ has
$\Gamma$-minimal elements; define $\cP_1$ to be the set of these
elements.  Let $\cN_1 = (\cP_1 + \Gamma) \cap \cB$.  Inductively,
having defined $\cP_1,\ldots,\cP_{n-1}$ and $\cN_1,\ldots,\cN_{n-1}$
for some $n \geq 2$, let $\cP_n$ consist of the $\Gamma$-minimal
elements of $\cB \minus \cP_{n-1}$, and set $\cN_n = (\cP_n + \Gamma)
\cap \cB$.
In other words,
\begin{itemize}
\item%
$\cP_n$ is the set of all positions $p \in \cB$ for which
$(p - \Gamma) \cap \cB$ is contained in $\cN_{n-1}$;
\item%
$\cN_n$ is the set of all positions $p \in \cB$ such that
$p - \gamma \in \cP_n$ for some $\gamma \in \Gamma$.
\end{itemize}
Note that $\cN_{n-2} \subseteq \cN_{n-1} \implies \cP_{n-1} \subseteq
\cP_n \implies \cN_{n-1} \subseteq \cN_n$, so it follows by induction
on~$n$, starting from $\cN_0 = \nothing$, that these containments all
hold.

\begin{lemma}\label{l:NP}
Let $\cP = \bigcup_{k=1}^\infty \cP_k$ and $\cN = \bigcup_{k=1}^\infty
\cN_k$.  Then
\begin{align*}
  \cP &= \{p \in \cB \mid (p - \Gamma) \cap \cB \subseteq \cN\}, \text{ and}
\\\cN &= \{p \in \cB \mid p - \gamma \in \cP \text{ for some } \gamma \in \Gamma\}.
\end{align*}
\end{lemma}
\begin{proof}
If $p \in \cP$, then $p \in \cP_k$ for some~$k$, so $(p - \Gamma) \cap
\cB \subseteq \cN_{k-1} \subseteq \cN$.  On the other hand, if $(p -
\Gamma) \cap \cB \subseteq \cN$ then $(p - \Gamma) \cap \cB \subseteq
\cN_{k-1}$ for some~$k$, since $\Gamma$ is finite, so $p \in \cP_k
\subseteq \cP$.

If $p \in \cN$ then $p \in \cN_k$ for some~$k$, so $p - \gamma \in
\cP_k \subseteq \cP$ for some $\gamma \in \Gamma$.  On the other hand,
if $p - \gamma \in \cP$ for some $\gamma \in \Gamma$, then $p - \gamma
\in \cP_k$ for some~$k$, so $p \in \cN_k \subseteq \cN$.
\end{proof}

For the existence claimed by the theorem, we check that the sets $\cP$
and $\cN$ from the lemma satisfy the axioms for sets of winning and
losing positions, respectively.
\begin{lemma}
With $\cP$ and~$\cN$ as in Lemma~\ref{l:NP}, we have $\cP \cap \cN =
\nothing$.
\end{lemma}
\begin{proof}
If $p \in \cP \cap \cN$ then $p \in \cP_n \cap \cN_{n-1}$ for
some~$n$, since the unions defining $\cP$ and~$\cN$ are increasing.
But $\cP_n \subseteq \cB \minus \cN_{n-1}$ by definition.
\end{proof}
\begin{lemma}
With $\cP$ and~$\cN$ as in Lemma~\ref{l:NP}, we have $\cP \cup \cN =
\cB$.
\end{lemma}
\begin{proof}
Suppose $p$ is $\Gamma$-minimal in $\cB \minus (\cP \cup \cN)$.  Then
$(p - \Gamma) \cap \cB \subseteq \cP \cup \cN$.  Therefore $p$ must
lie in~$\cP$ or in~$\cN$ by Lemma~\ref{l:NP}.
\end{proof}
By Lemma~\ref{l:NP}, it follows immediately that $(\cP + \Gamma) \cap
\cB = \cN$.

To prove uniqueness, suppose $\cB = \cP \cup \cN = \cP' \cup \cN'$,
where $\cP,\cN$ and $\cP',\cN'$ are pairs of winning and losing
positions.  First, suppose $\cP \cap \cN' = \nothing = \cP' \cap \cN$.
The first equality implies $\cN' \subseteq \cN$ while the second
equality implies $\cN \subseteq \cN'$.  Hence $\cN = \cN'$ and thus
$\cP = \cP'$.  Now suppose, by symmetry, that $\cP \cap \cN' \neq
\nothing$.  Let $p \in \cP \cap \cN'$.  If $p$ is not $\Gamma$-minimal
in~$\cB$, then since $p \in \cN'$, there is some $\gamma \in \Gamma$
such that $p - \gamma \in \cP'$; and since $p \in \cP$, we have $p -
\gamma \in \cN$ (note that $p - \gamma \in \cB$ since $\cP' \subseteq
\cB$).  Thus $q = p - \gamma \in \cP' \cap \cN$.  Continuing in this
manner and applying Lemma~\ref{l:dcc}, we reduce to the case where $p$
is $\Gamma$-minimal in~$\cB$.  But then $p \in \cP'$, contradicting $p
\in \cN'$.  Therefore $\cP \cap \cN' = \nothing$ and hence $\cP =
\cP'$ and $\cN = \cN'$.  This completes the proof.
\end{proof}

\section{Arbitrary impartial games as lattice games}\label{s:arbitrary}

As it turns out, lattice games can be rigged to encode arbitrary
games.  To make a precise statement, we briefly recall the definitions
of closed sets of games; see \cite{misereQuots} and its references for
more details.

Formally speaking, a finite impartial combinatorial game is a set
consisting of \emph{options}, each of which is, recursively, a finite
impartial combinatorial game.  The \emph{disjunctive sum} of two games
$G$ and $H$ is the game $G + H$ whose options are the union of $\{G' +
H \mid G'$ is an option of~$G\}$ and $\{G + H' \mid H'$ is an option
of~$H\}$.  A set of games is \emph{closed} if it is closed under
taking options and under disjunctive sum.  In particular, the
\emph{closure} of a single game~$G$ is the free commutative monoid
on~$G$ and its \emph{followers}, meaning the games obtained
recursively as an option, or an option of an option, and so on.  The
\emph{birthday} of a game is its longest chain of followers.

The theory of lattice games we develop here is, in the following
sense, universal.

\begin{thm}\label{t:arbitrary}
The closure of an arbitrary finite impartial combinatorial game, in
normal or mis\`ere play, can be encoded as a lattice game played
on~$\NN^d$.
\end{thm}
\begin{proof}
Let the game have $d$ distinct followers $G_2,\ldots,G_d$, and set
$G_1 = G$.  The disjunctive sum $G_{i_1} + \cdots + G_{i_r}$
corresponds to the position vector $e_{i_1} + \cdots + e_{i_r} \in
\NN^d$, where $e_1,\ldots,e_d$ are the standard basis vectors.  The
rule set $\Gamma$ consists of the vectors $e_i - e_j$ for each $i$
and~$j$ such that $G_j$ is an option of~$G_i$.  Normal play is encoded
by setting $\cD = \nothing$, and mis\`ere play is encoded by setting
$\cD = \{\0\}$.

It remains to verify the axioms for a rule set
(Definition~\ref{d:ruleset}).  The positivity axiom holds because the
function $\{1,\ldots,d\} \to \NN$ of the indices that sends $i$ to the
birthday of~$G_i$ is positive on~$\Gamma$.  The tangent cone axiom
holds by Lemma~\ref{l:path} because every position has a sequence of
moves to~$\0$.
\end{proof}

\begin{remark}
If any game can be encoded using~$\NN^d$, why allow arbitrary
polyhedral game boards?  Beyond Section~\ref{sub:polyhedral}, there
are at least two more answers.
\begin{enumerate}
\item%
For an arbitrary game, the encoding on $\NN^d$ via its game tree is
inefficient (exponentially so) in the sense of complexity theory;
polyhedral game boards allow efficient encodings of wider classes of
games.  (Even for a game on a board of the form~$\NN^d$ for some~$d$,
there is almost surely a better encoding as a lattice game than the
one provided in the proof of Theorem~\ref{t:arbitrary}.)
\item%
Our conjectures in Section~\ref{s:rational} concerning stratifications
in the polyhedral context attempt to place certain kinds of
periodicity in their natural generality, retaining only those
hypotheses we believe essential for the regularity to arise; there is
no reason, at present, to think that the simplicity of~$\NN^d$ has any
bearing.
\end{enumerate}
\end{remark}

\section{Squarefree games in normal play}\label{s:squarefree}

The notion of \emph{octal game} encompasses quite a broad range of
examples, but it can sometimes feel contrived, such as when
coincidences between the rules and certain heap sizes cause positions
with nonempty collections of heaps that nonetheless have no options
(see Example~\ref{e:octal}).  Lattice games suggest a natural common
generalization of octal games, as well as hexadecimal games and indeed
arbitrary heap games that, in particular, automatically avoids the
no-option phenomenon.  The \emph{squarefree games} we define here are
precisely those lattice games played on~$\NN^d$ such that the
Sprague--Grundy theorem for normal play finite impartial games is
commensurate with the coordinates placed on positions by virtue of the
game being on~$\NN^d$.

In this section, we assume that $\Pi = C = \RR_+^d$ and hence $\Lambda
= \NN^d$.  The following notation will come in handy a few times.

\begin{defn}
Given $v \in \RR^d$, let $v_+$ and $v_-$ be the nonnegative vectors
with disjoint support such that $v = v_+ - v_-$.  That is, $v_+ = v
\wedge \0$ and $v_- = -(v \vee \0) = v_+ - v$.
\end{defn}

\begin{prop}\label{p:squarefree}
For a rule set~$\Gamma$, the following are equivalent.
\begin{enumerate}
\item%
For each $\gamma \in \Gamma$ and $p \in \NN^d$, if $2p - \gamma \in
\NN^d$ then $p - \gamma \in \NN^d$.
\item%
If $p + p'$ is an option of $p + p$, then $p'$ is an option of $p$.
\item%
The maximum entry of each $\gamma \in \Gamma$ is equal to $1$.
\item%
The positive part $\gamma_+$ is a $0$-$1$ vector for all
$\gamma\in\Gamma$.
\item%
Each move takes away at most one heap of each size.
\end{enumerate}
\end{prop}
\begin{proof}
$1 \iff 2$: Assume 1 holds.  Suppose $p + p'$ is an option of $p + p$,
so $p' = p - \gamma$ for some $\gamma \in \Gamma$.  Then $p' \in
\NN^d$ and hence is an option of $p$.  Conversely, assume 2 holds.
Suppose $2p - \gamma \in \NN^d$.  Let $p' = p - \gamma$.  Then $p +
p'$ is an option of $p + p$, so $p'$ is an option of $p$, hence $p'
\in \NN^d$.

$1 \implies 3$: Suppose there exists $\gamma \in \Gamma$ with an entry
greater than 1.  Let $M = \text{max}\{\gamma_1,\ldots,\gamma_d\}$, and
let $p = \ceil{\frac M2}\1$ where $\1$ is the vector with all entries
equal to $1$.  Then the minimum of the entries of $2p - \gamma$ is $1$
for odd $M$ and $0$ for even $M$, and hence $2p - \gamma \in \NN^d$.
However, the minimum of the entries of $p - \gamma$ is $\ceil{\frac
M2} - M$ which is negative since $M > 1$.  Hence $p - \gamma \notin
\NN^d$.

$3 \implies 1$: Suppose 3 holds.  Let $\gamma \in \Gamma$ and let $p
\in \NN^d$ with $2p - \gamma \in \NN^d$.  For all $i$ such that $p_i =
0$, we must have $\gamma_i \leq 0$, and hence $(p - \gamma)_i \geq 0$.
For all $j$ such that $p_j > 0$, we have $(p - \gamma)_j = p_j -
\gamma_j \geq 1 - 1 = 0$.  Therefore $p - \gamma \in \NN^d$.

$3 \iff 4$: This is elementary, since every move $\gamma$ must possess
a strictly positive entry.

$4 \iff 5$: In the notation of Examples~\ref{e:nd} and~\ref{e:heaps},
condition~5 is the translation of condition~4 into the language of
heaps.
\end{proof}

\begin{defn}\label{d:squarefree}
We say a rule set $\Gamma$ is \emph{squarefree} if it satisfies any of
the equivalent conditions in Proposition~\ref{p:squarefree}.
\end{defn}

\begin{example}\label{e:octal}
The historically popular \emph{octal games}, invented by Guy and Smith
\cite{GuySmith56}, are heap games in which every move consists of
selecting a single heap and, depending on the heap's size and the
game's rules, either
\begin{enumerate}
\item%
removing the entire heap;
\item%
removing some beans from the heap, making it smaller; or
\item%
removing some beans from the heap, splitting it into two smaller heaps.
\end{enumerate}
A problem arises when there is a heap of size $k$ but the rules do not
allow the removal of $j$ beans from a heap, for any $j \leq k$.  In
this case, we may simply ignore heaps of size $k$ (treat them as heaps
of size $0$), and octal games naturally become a special class of
squarefree games.
\end{example}

\begin{defn}\label{d:normalplay}
The \emph{normal play} game board in~$\NN^d$ for a given rule set is
the one with no defeated positions: $\cD = \nothing$.
\end{defn}

Our results on normal play squarefree games are best stated in the
following terms.

\begin{defn}\label{d:cong}
Two positions $p,q \in \cB$ are \emph{congruent}, written $p \cong q$,
if
$$%
  (p + C) \cap \cP = p - q + (q + C) \cap \cP.
$$
In other words, $p + r \in \cP \iff q + r \in \cP$ for all $r$ in the
recession cone~$C$ of~$\cB$.
\end{defn}

\begin{excise}{%
  \begin{prop}\label{p:equiv}
  The relation $\cong$ as in Definition~\ref{d:cong} is an equivalence
  relation.
  \end{prop}
  \begin{proof}
  It is obvious that $\cong$ is reflexive.  Suppose $p \cong q$.
  Then, for any position $r$, we have $p + r \in \cP \iff q + r \in
  \cP$, so $q \cong p$.  Therefore $\cong$ is symmetric.  Finally,
  suppose $p \cong q$ and $q \cong r$.  Then
  \begin{eqnarray*}
  (p + \Lambda) \cap \cP
    &=& p - q + (q + \Lambda) \cap \cP \\
    &=& p - q + q - r + (r + \Lambda) \cap \cP \\
    &=& p - r + (r + \Lambda) \cap \cP.
  \end{eqnarray*}
  Hence $\cong$ is transitive.
  \end{proof}
  
  \begin{prop}\label{p:congadd}
  Let $\cB$ be a game board with winning positions~$\cP\!$, losing
  positions~$\cN$, and rule set $\Gamma$.  If $p \cong q$, then $p + r
  \cong q + r$ for any $r \in \Lambda$.
  \end{prop}
  \begin{proof}
  If $p \cong q$, then, in particular, $p + r + s \in \cP \iff q + r +
  s \in \cP$ for any $s \in \Lambda$.
  \end{proof}
  
  \begin{defn}\label{d:misereQuot}
  The \emph{quotient} of a lattice game played on a polyhedral cone
  $\Pi = C$ is the monoid $Q$ consisting of the equivalence classes
  from Proposition~\ref{p:equiv} along with the \emph{bipartition}
  of~$Q$ obtained by pushing down the winning and losing positions
  from~$\Lambda$.
  \end{defn}
}\end{excise}%

It is elementary to verify that congruence is an equivalence relation,
and that it is additive, in the sense that $p \cong q \implies p + r
\cong q + r$ for all $r \in C \cap \ZZ^d$.  Thus, when $\cB = \Lambda
= C \cap \ZZ^d$ is a monoid, the quotient of~$\cB$ modulo congruence
is again a monoid.

Throughout the remainder of this section, we assume that $\cB$ is the
normal play game board~$\NN^d$ with winning positions $\cP$, losing
positions~$\cN$, and rule set~$\Gamma$.

\begin{prop}\label{p:P=0}
If $p \in \cB$, then $p \in \cP \iff p \cong \0$.
\end{prop}
\begin{proof}
Suppose $p \in \cP$.  Let $q \in \cP$.  We claim that $p + q \in \cP$.
Clearly this is true for $p = \0$.  Now assume $p \succ \0$, and
suppose that $p' \in \cP \implies p' + q \in \cP$ for all $p' \prec
p$.  Let $\gamma \in \Gamma$ such that $p - \gamma \in \cB$.  Then $p
- \gamma \in \cN$, so there is $\gamma' \in \Gamma$ such that $p -
\gamma - \gamma' \in \cP$.  By our induction hypothesis, $(p + q -
\gamma) - \gamma' = (p - \gamma - \gamma') + q \in \cP$, so $p + q -
\gamma \in \cN$.  Since $\gamma$ was arbitrary, $p + q \in \cP$.
Hence $q \in \cP \implies p + q \in \cP$.

Now suppose $q \in \cN$.  Then there is $\gamma \in \Gamma$ such that
$q - \gamma \in \cP$.  Hence $p + q - \gamma \in \cP$, so $p + q \in
\cN$.  Therefore $p \cong \0$.

Conversely, suppose $p \cong \0$.  Then $p \in \cP$ since $\0\in\cP$.
\end{proof}

\begin{prop}\label{p:2p=0}
If\/ $\Gamma$ is squarefree and $p \in \cB$, then $2p \cong \0$.
\end{prop}
\begin{proof}
By Proposition \ref{p:P=0}, it suffices to show that $2p \in \cP$.  It
is clearly true for $p = \0$.  Now suppose $p \succ \0$ and $2p' \in
\cP$ for all $p' \prec p$.  Let $\gamma \in \Gamma$ such that $2p -
\gamma \in \cB$.  Since $\Gamma$ is squarefree, $p - \gamma \in \cB$,
hence $2(p - \gamma) \in \cB$.  By our induction hypothesis, $2p -
\gamma - \gamma = 2(p - \gamma) \in \cP$, hence $2p - \gamma \in \cN$.
Since $\gamma$ was arbitrary, $2p \in \cP$.
\end{proof}

\begin{cor}
If $p,q \in \cB$, then $p \cong q \iff p + q \in \cP$.
\end{cor}

\begin{lemma}\label{l:p=nq+r}
Let $n$ be a positive integer.  If $p \in \ZZ^d$, then there exist
unique $q \in \ZZ^d$ and $r \in \{0,\ldots,n-1\}^d$ such that $p = nq
+ r$.
\end{lemma}
\begin{proof}
Let $1 \leq i \leq d$.  There are unique $q_i \in \ZZ$ and $r_i \in
\{1,\ldots,n-1\}$ with $p_i = nq_i + r_i$ by the division algorithm.
Thus $q = (q_1, \ldots, q_d)$ and $r = (r_1, \ldots, r_d)$ do the job.
\end{proof}

\begin{thm}\label{t:normalplay}
Let $\Gamma$ be a squarefree rule set.  If $\cP_0 = \cP \cap
\{0,1\}^d$ then
$$%
  \cP = \cP_0 + 2\NN^d.
$$
\end{thm}
\begin{proof}
Let $w \in \cP_0$.  If $p \in 2\NN^d$, by Propositions~\ref{p:P=0}
and~\ref{p:2p=0}, $w + p \cong w \cong \0$ and hence $w + p \in \cP$.
On the other hand, let $w \in \cP$.  By Lemma~\ref{l:p=nq+r}, we may
write $w = 2p + q$ for some $p \in \cB$ and $q \in \{0,1\}^d$.  By
Propositions~\ref{p:P=0} and~\ref{p:2p=0}, we have $\0 \cong w \cong
2p + q \cong q$, hence $q \in \cP_0$.
\end{proof}

\begin{cor}\label{c:normalplay}
In normal play, if the rule set is squarefree, then the set of winning
positions is a finite disjoint union of translates of an affine
semigroup.
\end{cor}

\section{Algorithm for normal play}

In this section, we examine a method of computing the $\cP_0$ in
Theorem \ref{t:normalplay}.  To do this, we introduce a notion called
patterns.  For any set $S \subset \ZZ^d$, let $S_2$ denote the set of
elements of $S$ with entries modulo 2.

\begin{defn}
A \emph{pattern} is a subset of $\ZZ_2^d$, where $\ZZ_2 = \ZZ/2\ZZ$ is
the integers mod~$2$.
\end{defn}

\begin{defn}
Let $P$ be a pattern and let $N = \ZZ_2^d \minus P$.  We say that
$\Gamma$ \emph{sustains} $P$ if $P + \Gamma_2 = N$.
\end{defn}

\begin{prop}\label{p:sustain}
Let $\Gamma \subseteq \ZZ_2^d$.  If $\0 \notin \Gamma_2$, then
$\Gamma$ sustains a pattern $P$.  Furthermore, there is an algorithm
for finding $P$.
\end{prop}
\begin{proof}
Let $p_1 \in \ZZ_2^d$.  For each $k > 1$, having defined $p_j$ for $1
\leq j < k$, we choose $p_k \in \ZZ_2^d\,\minus\,\bigcup_{j=1}^{k-1}
(p_j + \Gamma_2)$.  Since $\ZZ_2^d$ is finite, this algorithm must
terminate after some $p_n$.  We claim that $\Gamma$ sustains $P =
\{p_1,\ldots,p_n\}$.  Let $N = \ZZ_2^d \minus P$.  Suppose $p_i = p_j
+ \gamma$ for some $i < j$ and some $\gamma \in \Gamma_2$.  Then $p_j
= p_i - \gamma = p_i + \gamma \in p_i + \Gamma_2$, a contradiction.
Hence $P + \Gamma_2 \subseteq N$.  Now suppose $p \in N$.  Then there
is some $k$ such that $p \in p_k + \Gamma_2 \subseteq P + \Gamma_2$.
\end{proof}

\begin{remark}
Note that the pattern~$P$ we obtain from the algorithm in
Proposition~\ref{p:sustain} is not necessarily unique, since it
depends on the choice of~$p_1$.
\end{remark}

\begin{thm}\label{t:alg}
If $\cB$ is a normal play game board with squarefree rule set~$\Gamma$
and winning positions $\cP = \cP_0 + 2\NN^d$, then there is an
algorithm for computing~$\cP_0$.
\end{thm}
\begin{proof}
Since $\Gamma$ is squarefree, $\0 \notin \Gamma_2$.  By the algorithm
in Proposition \ref{p:sustain}, we may obtain a pattern $P$ sustained
by $\Gamma$ such that $\0 \in P$.  We claim that $P = \cP_0$.  By
Theorem \ref{t:normalplay}, whether a position $p$ lies in $\cP$ or
$\cN$ depends solely on its coordinates modulo 2.  Therefore $p \in
\cP$ if and only if $p - \Gamma_2 \in \cN$.  Hence any pattern
sustained by $\Gamma$ is a viable candidate for $\cP_0$.  In
particular, $P$ works because it contains $\0$.  By Theorem
\ref{t:uniqueness}, $P$ is the only pattern that works.  Hence $\cP_0
= P$.
\end{proof}

\section{Rational strategies and affine stratifications}\label{s:rational}

As with any subset of $\NN^d$, the set $\cP$ of P-positions in a
lattice game can be recorded via its \emph{generating function}: the
power series $\sum_{p \in \cP} \ttt^p$ in the indeterminates $\ttt =
(t_1,\ldots,t_d)$,$\text{ where } \ttt^p = t_1^{\pi_1} \cdots
t_d^{\pi_d} \text{ for } p = (\pi_1,\ldots,\pi_d)$.  For particularly
well-behaved subsets of~$\NN^d$, the generating function is
\emph{rational}, meaning that it equals the Taylor series expansion of
a ratio of polynomials in $\ttt$ with integer coefficients.  These
notions make sense as well for generating functions supported on
pointed polyhedra, but the series are Laurent, not Taylor.

\begin{defn}
A \emph{rational strategy} for a lattice game is a rational generating
function for the set of P-positions.
\end{defn}

\begin{example}
Consider the game of \textsc{Nim} with heaps of size at most~$2$.
In normal play, a rational strategy is
$$%
  f(\cP;a,b) = \frac{1}{(1-a^2)(1-b^2)},
$$
the rational generating function for the affine semigroup~$2\NN^2$.
In mis\`ere play, a rational strategy is
$$%
  f(\cP;a,b) = \frac{a}{1-a^2} + \frac{b^2}{(1-a^2)(1-b^2)},
$$
with the first term enumerating every other lattice point on the first
axis, and the second enumerating the normal play P-positions that lie
off of the first axis.
\end{example}

\begin{example}\label{e:normalplay}
A squarefree game in normal play has a rational strategy
$$%
  f(\cP;\ttt) = \sum_{p\in\cP_0} \frac{\ttt^p}{(1-t_1^2)\cdots(1-t_d^2)},
$$
in the notation from Theorem~\ref{t:normalplay}.
\end{example}

Of what use is a rational strategy?  When one exists, it constitutes a
desirable data structure for representing and manipulating sets of
lattice points.

\begin{thm}[\cite{algsCGT}]\label{t:rational}
Any rational strategy for a lattice game, presented as a ratio of two
polynomials with integer coefficients, produces algorithms for
\begin{itemize}
\item%
determining whether a position is P-position or an N-position, and
\item%
computing a legal move to a P-position, given any N-position.
\end{itemize}
These algorithms are efficient when the rational strategy is a
\emph{short} rational function, in the sense of Barvinok and Woods
\cite{BaWo03}.
\end{thm}

Theorem~\ref{t:rational} is proved by straightforward application of
the algorithms of Barvinok and Woods \cite{BaWo03}.  That said, the
subtlety lies more in phrasing the statement precisely, particularly
when it comes to complexity; see \cite{algsCGT}.

Lattice games need not a~priori possess rational strategies, but
examples and heuristic arguments lead to the following assertion.

\begin{conj}\label{conj:rational}
Every lattice game possesses a rational strategy.
\end{conj}

Conjecture~\ref{conj:rational} is precise, but it allows for such a
wide array of generating functions that it fails to capture the
regularities pervading all examples to date.  In fact, we were led to
Conjecture~\ref{conj:rational} only after observing the existence of
certain decompositions.

\begin{defn}\label{d:strat}
An \emph{affine stratification} of a lattice game is a partition
$$%
  \cP = \biguplus_{i=1}^r W_i
$$
of the set of P-positions into a disjoint union of sets $W_i$, each of
which is a finitely generated module for an affine semigroup $A_i
\subset \ZZ^d$; that is, $W_i = F_i + A_i$, where $F_i \subset \ZZ^d$
is a finite set.
\end{defn}

\begin{example}
Consider again the game of \textsc{Nim} with heaps of size at
most~$2$.  An affine stratification for this game is $\cP = 2\NN^2$;
that is, the entire set of P-positions forms an affine semigroup.  In
mis\`ere play, $\cP = \big((1,0) + 2\NN\big) \uplus \big((0,2) +
2\NN^2\big)$ is the disjoint union of $W_1 = 1 + 2\NN$ (along the
first axis) and $W_2$, which equals the translate by twice the second
basis vector of the affine semigroup~$2\NN^2$.
\end{example}

\begin{example}
Every normal play squarefree game has an affine stratification; this
is Corollary~\ref{c:normalplay}.
\end{example}

\begin{conj}\label{conj:strat}
Every lattice game possesses an affine stratification.
\end{conj}

\begin{remark}
Conjecture~\ref{conj:strat} is equivalent to the same statement with
the extra hypothesis that the rule set is saturated.  Indeed,
$\ZZ\Gamma$ has finite index in~$\ZZ^d$, whence the game board is a
disjoint union of finitely many games each of whose rule sets is
saturated in its ambient lattice.
\end{remark}

The importance of Conjecture~\ref{conj:strat} stems from its
computational consequences.

\begin{thm}[\cite{algsCGT}]\label{t:strat}
A rational strategy can be efficiently computed from any affine
stratification.
\end{thm}

Again, the proof comes down to the algorithms of Barvinok and Woods
\cite{BaWo03}, but the notion of ``efficiency'' must be made precise,
and that is even more subtle than Theorem~\ref{t:rational}; see
\cite{algsCGT}.

\begin{example}
The mis\`ere lattice game on~$\NN^5$ whose rule set forms the
columns~of
$$%
\newcommand\ph{\phantom{-}}
\Gamma = \left[
         \begin{array}{@{}*{8}{r@{\ \,}}@{}r}
           \!\ph1&\ph0&\ph0&\ph0&-1& 0& 0& 0\\
           \!   0&   1&   0&   0& 1&-1& 0& 0\\
           \!   0&   0&   0&   0& 0& 1&-1& 0\\
           \!   0&   0&   1&   0& 0& 0& 1&-1\\
           \!   0&   0&   0&   1& 0& 0& 0& 1
         \end{array}
         \right]
$$
has infinite mis\`ere quotient \cite[Section~A.7]{misereQuots}.  The
illustration of the winning positions in this lattice game provided by
Plambeck and Siegel \cite[Figure~12]{misereQuots} was one of the
motivations for the definitions in this paper, particularly
Definition~\ref{d:strat}, because it possesses an interesting
description as an affine stratification.  Indeed, for this lattice
game, $\cP = W_1 \uplus \cdots \uplus W_7$ for modules $W_k = F_k +
A_k$ over the afffine semigroups
\begin{align*}
A_1 &= \NN\{(2,0,0,0,0),(0,2,0,0,0),(0,0,1,0,0),(0,0,0,2,0),(0,0,0,0,2)\}\\
A_2 &= \NN\{(2,0,0,0,0),(0,2,0,0,0),(0,0,0,2,0),(0,0,0,0,2)\}\\
A_3 &= \NN\{(2,0,0,0,0),(0,2,0,0,0),(0,0,0,0,2),(0,0,0,2,2)\}\\
A_4 &= \NN\{(2,0,0,0,0),(0,2,0,0,0),(0,0,0,2,0),(0,0,0,2,2)\}\\
A_5 &= \NN\{(2,0,0,0,0),(0,2,0,0,0)\}\\
A_6 &= \NN\{(2,0,0,0,0),(0,2,0,0,0),(0,0,0,2,2)\}\\
A_7 &= \NN\{(2,0,0,0,0)\},
\end{align*}
where the finite generating sets $F_k$ consist of the columns of the
following:
\begin{align*}
&
  F_1 \colon \left[
        \begin{array}{c}
          \!0\!\\
          \!0\!\\
          \!3\!\\
          \!0\!\\
          \!0\!
        \end{array}
        \right]
  \quad
  F_2 \colon \left[
        \begin{array}{*{5}{c@{\ \,}}@{}c}
          \!0&1&1&1&1\!\\
          \!0&1&1&1&1\!\\
          \!2&2&2&2&2\!\\
          \!0&3&3&3&1\!\\
          \!0&1&3&5&5\!
        \end{array}
        \right]
  \quad
  F_3 \colon \left[
        \begin{array}{*{10}{c@{\ \,}}@{}c}
          \!1&0 &0&1&0&1&0 &1&0 &1 \!\!\\
          \!0&1 &0&0&1&1&0 &0&1 &1 \!\!\\
          \!2&2 &1&1&1&1&0 &0&0 &0 \!\!\\
          \!0&1 &0&0&1&1&0 &0&1 &1 \!\!\\
          \!9&12&9&8&9&8&10&9&12&13\!\!
        \end{array}
        \right]
\\
&
  F_4 \colon \left[
        \begin{array}{*{32}{c@{\ \,}}@{}c}
          \!1&1&1&1&0&0&0&0&0&0&1&1&1&0&0&0&1&1&0&0&0&0&0&1&1&1&1&0&0&0&1&1&1\!\\
          \!0&0&0&0&1&1&1&0&0&0&0&0&0&1&1&1&1&1&0&0&0&0&0&0&0&0&0&1&1&1&1&1&1\!\\
          \!2&2&2&2&2&2&2&1&1&1&1&1&1&1&1&1&1&1&0&0&0&0&0&0&0&0&0&0&0&0&0&0&0\!\\
          \!1&1&1&1&0&0&0&0&2&0&1&3&1&2&2&2&3&1&2&4&4&4&4&3&5&5&3&2&4&2&3&3&3\!\\
          \!0&2&4&6&1&5&7&0&4&6&0&4&6&1&5&7&1&5&0&2&4&6&8&0&2&4&6&1&5&7&1&5&7\!
        \end{array}
        \right]
\\
&
  F_5 \colon \left[
        \begin{array}{*{29}{c@{\ \,}}@{}c}
          \!1&0&1&1&1&1&0&1&1&1&1&0&0&0&0&0&0&1&1&1&1&1&1&0&0&0&0&0&1&1\!\\
          \!1&0&0&0&0&0&1&1&1&1&1&2&0&0&0&0&0&2&0&0&0&0&0&1&1&1&1&1&3&1\!\\
          \!2&1&1&1&1&1&1&1&1&1&1&0&0&0&0&0&0&0&0&0&0&0&0&0&0&0&0&0&0&0\!\\
          \!0&0&0&1&0&1&1&1&2&0&0&0&1&2&0&1&2&0&2&3&0&1&2&1&0&1&0&1&1&2\!\\
          \!3&3&2&3&5&6&4&1&3&4&7&0&2&3&4&5&6&1&2&3&3&4&5&0&2&3&5&6&1&4\!
        \end{array}
        \right]
\\
&
  F_6 \colon \left[
        \begin{array}{*{7}{c@{\ \,}}@{}c}
          \!0&0&1&1&0&0&1 \!\!\\
          \!0&0&0&0&1&1&1 \!\!\\
          \!0&0&0&0&0&0&0 \!\!\\
          \!0&1&0&1&0&1&1 \!\!\\
          \!7&8&6&7&8&9&10\!\!
        \end{array}
        \right]
  \quad
  F_7 \colon \left[
        \begin{array}{*{3}{c@{\ \,}}@{}c}
          \!0&1&1\!\\
          \!0&0&0\!\\
          \!0&0&0\!\\
          \!0&0&1\!\\
          \!1&0&1\!
        \end{array}
        \right].
\end{align*}
We are hopeful that this affine stratification, as a mode for
presenting the P-positions, will lead to an algorithmic verification
of the presentation for the mis\`ere quotient.
\end{example}


\end{document}